\documentclass[letterpaper, 10pt, conference]{ieeeconf}
\usepackage{graphicx, amsmath, bm}
\usepackage{srcltx}
\newtheorem{assumption}{Assumption}
\newtheorem{lemma}{Lemma}
\IEEEoverridecommandlockouts                
\title{\LARGE \bf
An Extended Kalman Filter with a Computed Mean Square Error Bound}
\author{G. Hexner, and  H. Weiss$^*$ 
\thanks{$^*$ The authors are with RAFAEL, Advanced Defense Systems,
{\tt\small \{georgeh, haimw\}@rafael.co.il}}%
}
\begin{document}
\maketitle
\thispagestyle{empty}
\pagestyle{empty}
\begin{abstract}
The paper proposes a new recursive filter for non-linear systems that
inherently computes a valid bound on the mean square estimation error.
The proposed filter, 
bound based extended Kalman,
(BEKF) is in the form of an extended Kalman filter.
The main difference of the proposed
filter from the conventional extended Kalman filter is in the use of
a computed mean square error bound matrix, to calculate the filter
gain, and 
to serve as bound on the actual mean square error. The paper shows
that when the system is linear the proposed filtering algorithm
reduces to the conventional Kalman filter.
The
theory presented in the paper is applicable to a wide class of
systems, but if the system is polynomial, then the recently 
developed theory of positive polynomials considerably simplifies
the filter's implementation.
\end{abstract}
\section{Introduction}
The purpose of this paper is to propose a new form of extended Kalman
filter (BEKF) with a computed bound on the mean squared error matrix,
used to calculate the 
filter gain and to serve as a bound on the actual estimate mean square
error. The
theory presented is general, subject only to some asymptotic growth
and continuity constraints, but  implementation is very much
simplified if the underlying system consists of rational polynomials.
In this case the recently introduced theory of positive polynomials
\cite{PapaPrajna}, \cite{ParriloLall}, \cite{Chesi}, and the software
SOSTOOLS \cite{SOSTOOLS} provides the tools for efficient implmentation
of the filter.
SOSTOOLS  translates the problem to a
semi-definite program, which is readily solved by
SeDuMi \cite{SeDuMi}. This set of software makes possible the
numerical calculation of the bounds necessary in this paper
routine.

Previous attempts at extended Kalman filtering for polynomial based
systems inevitably faced the closure problem \cite{Sorenson},
\cite{Jazwinsky}. The closure problem refers to the fact that to
calculate the $n$th moment of a distribution, the value of the $n+1$ and
possibly higher order moments are required. A
popular method has been to assume that moments of higher order are
related to lower order moments as if the underlying probability density
were Gaussian, \cite{Jazwinsky}.

One approach to non-linear estimation was presented in
\cite{KumarChakra} based on a special type of  discretization of the
exact equations of nonlinear filtering. 
A different approach to estimation for non-linear systems was proposed
for cone bounded non-linearities in
\cite{GilmanRhodes1,GilmanRhodes2,GilmanRhodes3}. The 
special feature of these papers, compared to the many publications
that deal with estimation for non-linear systems, is the derivation of
an {\em analytic bound} on the performance of the estimator, without
requiring any sort of truncation approximation. The present paper
proposes to derive an analytic bound on the mean square
estimation error using an alternative approach.

For proper operation of any Kalman filter it is essential that the
calculated filter mean square error track the actual filter mean
square error reasonably well. The reason for this is that the filter
mean square error defines the filter gain. Too small  value of the
calculated mean square error implies that the Kalman gain is too low
and the observations are insufficiently weighted in updating the
filter estimate. The formalization of the concept is called
consistency, \cite{BarShalom}. For a filter to be consistent two
conditions have to be fulfilled:
\begin{enumerate}
\item Have mean zero (i. e. the estimates are unbiased)
\item Have covariance matrix as calculated by the filter.
\end{enumerate}
The exact methods for testing and consequences of filter consistency
are discussed 
further in \cite{BarShalom}. The proposed filter (BEKF) does not
ensure that the estimate error is zero mean, but does ensure 
that the filter mean square error calculated by the algorithm
dominates the actual mean square error. Thus it cannot be claimed that
BEKF is consistent; however, it does ensure a
reasonable value for the filter gain.

An important step in the development of an extended Kalman filter is
``tuning''. This consists of adjusting (usually increasing) by trial
and error the
intensity of the process noise and possibly the observation noise so
that the filter calculated mean square error is in some agreement with
the actual mean square error. The contribution of the present paper is
the development of an algorithm that ensures that the calculated mean
square estimation error is larger or equal to the actual estimation
error. In particular the new filter precludes the possibility filter
divergence. 
\section{Problem Statement}
Let $\bm{x}$ be an $N$ dimensional diffusion,
\begin{equation}\label{eq:dif}
{\rm d}\bm{x} = \bm{f(x)}{\rm d}t + \bm{g(x)}{\rm d}\bm{w}
\end{equation}
where $\bm{w}$ is a standard vector Weiner process. The variable $\bm{y}$
is observed  at time instances $T_k,\:k=0,1,\ldots$ 
\begin{equation}\label{eq:obseq}
\bm{y}(T_k) = \bm{Hx}(T_k) + \bm{v}_k
\end{equation}
where $\bm{v}_k$ is a sequence of zero mean independent  variables
with covariance matrix  $\bm{R}$. Also the initial mean squared error
matrix $\bm{\Sigma}(0)$ of the initial value of the state vector
$\bm{x}(0)$ is assumed known.
The following is necessary for the development of the bound:
\begin{assumption}\label{as:ContDif}
The function $\bm{g(x)}$ is continuous and the function $\bm{f(x)}$ is
continuously differentiable. 
\end{assumption}
\begin{assumption}\label{as:Pieq}
For any symmetric matrix $\bm{P}$ there exists a symmetric matrix
$\bm{Q}$ and a constant $q$ such that 
\begin{equation}\label{eq:assumeq}
\left(\frac{\partial \bm{f}(\bm{x})}{\partial \bm{x}}\right)^\prime
\bm{P} + \bm{P}\left(\frac{\partial \bm{f}(\bm{x})}{\partial \bm{x}}\right)
\preceq \bm{Q} 
\end{equation}
and
\begin{equation}\label{eq:assumeq2}
{\rm Tr}\{\bm{g^\prime(x) P g(x)}\} \le q
\end{equation}
\end{assumption}
Note that   there are no sign definite constraints on $\bm{P}$ or
$\bm{Q}$ in (\ref{eq:assumeq}). Assumption \ref{as:Pieq} essentially
limits the growth of $\|\bm{f(x)}\|$ to be at most linear in $\bm{x}$,
for large $\|\bm{x}\|$.
The form on the left hand side of (\ref{eq:assumeq}) occurs in the
study of contracting systems \cite{Lohmiller}. When $\bm{P} \succ
\bm{0}$ the form measures the rate that two solutions of
(\ref{eq:dif}), with close initial conditions diverge from each
other. In \cite{Pham} it is shown that if $\bm{P} \succ \bm{0}$ and if 
$\bm{Q} \prec \bm{0}$ in
(\ref{eq:assumeq}) in the whole space and (\ref{eq:assumeq2}) is valid
then (\ref{eq:dif}) is incrementally stochastically stable.

The aim of the paper is to develop an algorithm for calculating an
estimate, $\bm{\hat x}$ of the state of (\ref{eq:dif}), based  on the
observations (\ref{eq:obseq})  and a guaranteed bound \footnote{Note
  that the emphasis here is the development of a bound. No attempt is
  made here to develop tightest possible bound} for the expected
value of the mean square error of the estimate. 
\subsection{Some Examples of Systems Satisfying Assumption \ref{as:Pieq}}
 The simplest system satisfying Assumption \ref{as:Pieq} is
\begin{equation}
{\rm d}\bm{x} = \bm{A(x)x}{\rm d}t + \bm{g(x)}{\rm d}\bm{w}
\end{equation}
where $\bm{A(x)}$ and $\bm{g(x)}$ satisfy
\begin{equation}
||\bm{A(x)}|| + \left \|\frac{\partial
  \bm{A}(\bm{x})}{\partial \bm{x}}\right\|  \le M_{\bm{A}};\;
\;||\bm{g(x)}|| \le M_{\bm{g}} 
\end{equation}
for some constants $M_{\bm{A}}$, and $M_{\bm{g}}$. 

A system which contradicts assumption \ref{as:Pieq} is
\begin{equation}\label{eq:NotAssmp1}
{\rm d}x = (1 + x^2) dt + dw
\end{equation}
Note that the system (\ref{eq:NotAssmp1}) has a finite escape time.

\section{Time Update}
Between observations the state estimate, $\bm{\tilde x}$ evolves as in
the extended Kalman filter, according to 
\begin{equation}\label{eq:Tupdt}
{\rm d}\bm{\tilde x = f(\tilde x)}{\rm d}t
\end{equation}
The initial condition for $\bm{\tilde x}$ after the $k$th data processing
step is $\bm{\hat x}_k$.
Then the estimate error, $\bm{\tilde e = \tilde x - x}$ evolves as
\begin{equation}
{\rm d}\bm{\tilde e} = \bm{F(\tilde e, \tilde x)}{\rm d}t +
\bm{G(\tilde e, \tilde x)}{\rm d}\bm{w} 
\end{equation}
where
\begin{equation}
\bm{F(\tilde e, \tilde x)} = \bm{f(\tilde x) - f(\tilde x - \tilde e)}
\end{equation}
\begin{equation}
\bm{G(\tilde e, \tilde x)} = -\bm{g(\tilde x - \tilde e)}
\end{equation}
Note that $\bm{\tilde x}$ is computed according to (\ref{eq:Tupdt}),
so that it is a known function. 

An important contribution of the paper is the calculation of $\tilde
\Sigma$ such that 
\begin{equation}
{\rm \bf E}\{\bm{\tilde e \tilde e^\prime}\} = \bm{\Sigma} \preceq
\bm{\tilde \Sigma} 
\end{equation}
is a bound on the mean square error, ${\rm \bf E}\{\bm{\tilde e \tilde
  e}^\prime\} = \bm{\Sigma}$. 
This is accomplished by calculating bounds for
\begin{equation}
{\rm Tr}\{\bm{P \dot \Sigma}\}
\end{equation}
Using the Ito calculus \cite{Jazwinsky},
\begin{IEEEeqnarray}{rCl}
{\rm d}\left[\bm{\tilde e^\prime P \tilde e}\right]&  = & {\rm
  d}\bm{\tilde e^\prime P \tilde e + \tilde e^\prime P}{\rm
  d}\bm{\tilde e} + {\rm d} \bm{\tilde e^\prime} \bm{P} {\rm
  d}\bm{\tilde e} \\ \label{eq:dxPx} 
&= & \left[\bm{F^\prime(\tilde e, \tilde x) P \tilde e + \tilde
    e^\prime P F(\tilde e, \tilde x)}\right] \nonumber \\  
& &+\: {\rm Tr}\{\bm{G^\prime(\tilde e, \tilde x) PG(\tilde e, \tilde
  x)}\} {\rm d}t  \\  
& &+ \: {\rm d}\bm{w}^\prime \bm{G^\prime(\tilde e, \tilde x) P \tilde
  e + \tilde e^\prime 
    P G(\tilde e, \tilde x) }{\rm d}\bm{w}  \nonumber 
\end{IEEEeqnarray}
Taking expectation in (\ref{eq:dxPx}), and simplifying 
\begin{IEEEeqnarray}{rCl}
{\rm Tr}\{\bm{P \dot \Sigma}\}& =& {\rm \bf E}\{\bm{F^\prime(\tilde e, \tilde x) P \tilde e + \tilde e^\prime P F(\tilde e, \tilde x)} \nonumber \\ \label{eq:Edot}
& &\:+ {\rm Tr}\{\bm{G^\prime(\tilde e, \tilde x) P G(\tilde e, \tilde x)}\}\}
\end{IEEEeqnarray}
Using the mean value theorem,
\begin{equation}\label{eq:mvth}
\bm{F(\tilde e, \tilde x)}  =  \bm{f(\tilde x) - f(\tilde x - \tilde e)} 
 = \bm{\frac{\partial f(c)}{\partial x}\tilde e}
\end{equation}
where $\bm{c}$ is a point on the line connecting $\bm{\tilde x}$ and
$\bm{\tilde x - \tilde e}$. 
Substituting (\ref{eq:mvth}) in (\ref{eq:Edot}),
\begin{equation}
\begin{split}
{\rm Tr}\{\bm{P \dot \Sigma}\} =  {\rm \bf E}\{\bm{\tilde
  e^\prime\left(\frac{\partial  f}{\partial x}\right)^\prime P \tilde
  e + \tilde e^\prime P \left(\frac{\partial  f}{\partial x}\right)
  \tilde e}\} \\ 
+ {\rm Tr}{\rm \bf E}\{\bm{G^\prime(\tilde e, \tilde x) P G(\tilde e, \tilde x)}\}
\end{split}
\end{equation}
taking expectations and using (\ref{eq:assumeq}) yields the inequality
\begin{equation}\label{eq:Trineq}
{\rm Tr}\{\bm{P \dot \Sigma}\} \le {\rm Tr}\{\bm{Q \Sigma}\} + q
\end{equation}
Therefore assumptions \ref{as:ContDif} and \ref{as:Pieq} ensure that
for every $\bm{P}$ there exist a $\bm{Q}$ and $q$ satisfying
(\ref{eq:Trineq}).

In the next subsection a bound $\bm{\dot{\overline \Sigma}}$ for
$\bm{\dot \Sigma}$ is calculated based on the repeated use of
(\ref{eq:Trineq}), using a set of $\bm{P}_i$, $\bm{Q}_i$, and $q_i$. 
Then to each value of $\bm{\Sigma}$ the algorithm to be presented in
the next subsection calculates a bound $\bm{\dot{\overline \Sigma}}$
for $\bm{\dot{\Sigma}}$. That is, the algorithm defines a function,
$\bm{\dot{\overline\Sigma}}(\bm{\Sigma}, t)$ .
In
the succeeding subsection, based on this function, a differential
equation  
for the bound $\bm{\tilde \Sigma}$ for $\bm{\Sigma}$ is derived.

\subsection{The Bound for $\bm{\dot \Sigma}$}
Given a value for the
mean squared error matrix, $\bm{\Sigma}$ the procedure, to be presented
in this section,  calculates a bound for its 
derivative. Here, the bound for the derivative is denoted as
$\bm{\dot{\overline \Sigma}}$. Although, the bound depends on
$\bm{\Sigma}$ and $t$, and hence is a function of $\bm{\Sigma}$ and $t$,
in the
present subsection, this dependence is suppressed.
The calculation of the bound is carried out in
two steps:
\begin{enumerate}
\item A bounded set ${\cal S}$ containing $\bm{\dot \Sigma}$ is
  calculated.
\item A single $\bm{\dot {\overline \Sigma}}$ such that
\begin{equation}\label{eq:dotbarSigdef}
\bm{\dot \Sigma} \preceq \bm{\dot{ \overline \Sigma}}, \forall \bm{\dot
  \Sigma} \in {\cal S}
\end{equation}
\end{enumerate}
is calculated.
\subsubsection{Calculation of ${\cal S}$}
The set ${\cal S}$ is defined as
\begin{equation}\label{eq:SetS}
{\cal S}(\bm{\Sigma}) = \{\bm{\dot \Sigma}|{\rm Tr}\{\bm{P}_i\bm{\dot
  \Sigma}\} \le 
{\rm Tr}\{\bm{Q}_i\bm{\Sigma}\}+q_i, \forall i\} 
\end{equation}
where the set of $N\times N$ symmetric matrices $\bm{P}_i$ satisfy the
following condition: for any symmetric 
$N\times N$, non-zero  $\bm{Z}$ there exists an $i$ such that
\begin{equation}\label{eq:Pidef}
{\rm Tr}\{\bm{P}_i \bm{Z}\} > 0
\end{equation}
This condition ensures that if $\bm{\dot \Sigma}$ becomes unbounded,
 at least one of the inequalities in (\ref{eq:SetS}) is violated. 
Note that ${\cal S}$ is a  polytope. The
construction of a set of
$\bm{P}_i$ satisfying (\ref{eq:Pidef}) is presented in appendix
\ref{sec:PiCalc}, and the algorithm for calculating the corresponding
$\bm{Q}_i$, and $q_i$ is presented in section \ref{eq:QqCalc}.
\subsubsection{Calculation of $\bm{\dot{\overline \Sigma}}$}\label{sec:dotSigma}
This section relies on generalized inequalities with respect various
proper cones, see \cite{Boyd}. The generalized inequality  with
respect to the cone 
${\cal A}$ is denoted as $\preceq_{\cal A}$. Ideally a matrix
$\bm{\dot{\overline \Sigma}}$ is sought such that
\begin{equation}\label{eq:dotbarSigdef2}
\bm{\dot \Sigma} \preceq_{\cal M_+} \bm{\dot{\overline \Sigma}}, \;
\forall \bm{\dot \Sigma} \in {\cal S}
\end{equation}
where ${\cal M_+}$ is the cone of positive semi-definite matrices.
(This is the same as (\ref{eq:dotbarSigdef}), while explicitly
indicating the cone 
for the generalized inequality.) This is an infinite dimensional
problem. Let ${\cal F_+}$ denote the cone with a {\em finite} set of
generators, $\bm{U}_i\succeq_{\cal M_+} \bm{0}$. Then any 
$\bm{\dot{\overline \Sigma}}$ that satisfies
\begin{equation}
\bm{\dot \Sigma} \preceq_{\cal F_+} \bm{\dot{\overline \Sigma}}, \;
\forall \bm{\dot \Sigma} \in {\cal S}
\end{equation}
also satisfies (\ref{eq:dotbarSigdef2}), since 
${\cal F_+} \subseteq{\cal M_+}$. 
Any set of $\bm{U}_i\succeq_{\cal M_+} \bm{0}$ yields a
valid bound, $\bm{\dot{\overline \Sigma}}$, but the bound becomes
tighter as ${\cal F_+}$ approaches ${\cal M_+}$.

From the geometry of ${\cal M_+}$ for maximal ${\cal F_+}$, the
$\bm{U}_i$ should be on the boundary of ${\cal M_+}$, that is, the
$\bm{U}_i$ should be rank one matrices. For the example for second
order systems a reasonable choice is the four positive semi-definite
rank one matrices 
\begin{equation}
\left[\begin{array}{cc}
1 & 0\\
0 & 0
\end{array}\right], 
\left[\begin{array}{cc}
0 & 0\\
0 & 1
\end{array}\right], 
\left[\begin{array}{cc}
1 & 1\\
1 & 1
\end{array}\right], 
\left[\begin{array}{cc}
1 & -1\\
-1 & 1
\end{array}\right]
\end{equation}
These matrices served as the $\bm{U}_i$ for the example in section
\ref{sec:example}.
Note that there exist positive definite matrices that cannot be
expressed as a positive combination of these $\bm{U}_i$, one example
of such a positive definite matrix is
\begin{equation}
\left[\begin{array}{cc}
1 + 2\epsilon & 1 \\
1             & \frac{1+\epsilon}{1 + 2\epsilon}
\end{array}\right]
\end{equation}
for any $\epsilon > 0$. That is ${\cal F_+} \subset {\cal M_+}$.

Let $\bm{T}_i$ be the finite set of generators for the dual cone
${\cal F^*_+}$. The calculation of the $\bm{T}_i$,
from a given set of $\bm{U}_i$ is discussed in
appendix \ref{sec:TiCalc}.
From the properties of dual cones, \cite{Boyd} since 
${\cal F_+} \subseteq {\cal M_+}$ then
${\cal F^*_+} \supseteq {\cal M_+^*}={\cal M_+}$. Therefore,
for any $\bm{T} \succeq \bm{0}$ there exist $\lambda_i \ge 0$ such that
\begin{equation}\label{eq:anyT}
\bm{T} = \sum_i \lambda_i \bm{T}_i
\end{equation}
 Define $t_i$ 
\begin{equation}\label{eq:t_idef}
t_i = \max_{\bm{\dot \Sigma} \in {\cal S}(\bm{\Sigma})} {\rm
  Tr}\{\bm{T}_i\bm{\dot \Sigma}\} 
\end{equation}
Since the set ${\cal S}$ is bounded, 
the $t_i$ are finite. Define $\bm{\dot{\overline \Sigma}}$
\begin{equation}\label{eq:bardotSig}
\bm{\dot {\overline \Sigma}} = \sum_i s_i \bm{T}_i
\end{equation}
for some $s_i$ yet to be specified, and define the matrix $\bm{L}$
with entries $\bm{L}_{ij}$,
\begin{equation}\label{eq:Lcalc}
\bm{L}_{ij} = {\rm Tr}\{\bm{T}_i \bm{T}_j\}
\end{equation}
The matrix $\bm{L}$ is square with dimension equal to the number of
$\bm{T}_i$'s.  
Define the vectors $\bm{s}$ and $\bm{t}$, whose entries are $s_i$ and
$t_j$, respectively. Then to any $\bm{s}$ that satisfies 
\begin{equation}\label{eq:sdef}
\bm{t} \preceq \bm{Ls}
\end{equation}
there corresponds a bound $\dot{\bm{\overline \Sigma}}$ from the
assumed form (\ref{eq:bardotSig}). The calculation of  $\bm{s}$ is
discussed at the end of the next subsection.

That the $\bm{\dot{ \overline \Sigma}}$ defined in
(\ref{eq:bardotSig}) satisfies (\ref{eq:dotbarSigdef}) follows if it
can be shown that for any $\bm{T} \succeq \bm{0}$ 
\begin{equation}
{\rm Tr}\{\bm{T(\dot{\Sigma} - \dot{\overline \Sigma})}\} \le 0
\end{equation}
Indeed, from (\ref{eq:anyT}) and (\ref{eq:sdef}),
\begin{IEEEeqnarray}{rCl}
& &{\rm Tr}\left\{\bm{T(\dot{\Sigma} - \dot{\overline \Sigma})}\right\}\\
& = & {\rm Tr}\left\{\sum_i \lambda_i \bm{T}_i\left(\dot{\bm{\Sigma}}-\dot{\bm{\overline \Sigma}}\right)\right\}\\
& \le &\sum_i\lambda_i \left(t_i - \sum_j\bm{L}_{ij} \bm{s}_j\right) \\
&\le& 0
\end{IEEEeqnarray}
\subsection{The Mean Square Bound}
For each $\bm{\Sigma}$ a  bound $\bm{\dot{\overline \Sigma}}$ for
$\bm{\dot \Sigma}$ was
derived in the previous section. The purpose of the present section is
to calculate $\bm{\tilde \Sigma}$ such that
\begin{equation}
\bm{ \Sigma}(t) \preceq \bm{\tilde\Sigma}(t), \; \forall t
\end{equation}
The $\bm{\dot{\overline \Sigma}}$, calculated in the previous section
depends on $\bm{\Sigma}$, and through $\bm{\tilde x}$ on $t$. In the
present section this dependence is 
made explicit, and the bound is now denoted as 
$\bm{\dot{\overline    \Sigma}(\Sigma,t)}$.  
If $\bm{\Sigma}$ were available then solving the differential equation
\begin{equation}
\bm{\dot{\tilde \Sigma}} = \bm{\dot{\overline \Sigma}}(\bm{\Sigma}, t)
\end{equation}
would yield a bound $\bm{\tilde \Sigma}$. 
The differential
equation obtained by replacing $\bm{\Sigma}$ by $\bm{\tilde \Sigma}$
in $\bm{\dot{\overline\Sigma}( \Sigma, t)}$, yields an implementable
differential equation,
\begin{equation}\label{eq:DtildeSig}
\bm{\dot{\tilde \Sigma}} = \bm{\dot{\overline \Sigma}(\tilde \Sigma, t)}
\end{equation}
Given an initial condition it has a well defined solution.
The following lemma is an
extension of the comparison 
principle \cite{Khalil}, p. 102 to matrix valued functions, and is
used here to show that the solution of (\ref{eq:DtildeSig}) provides a
valid bound for $\bm{\Sigma}$.
\begin{lemma}\label{lem:comp}
Consider the matrix valued differential equation (\ref{eq:DtildeSig})
with initial condition
\begin{equation}
\bm{ \Sigma}(T_i) \preceq \bm{\tilde\Sigma}(T_i)
\end{equation}
Then
\begin{equation}
\bm{\Sigma}(t) \preceq \bm{\tilde \Sigma}(t), \; t \ge T_i
\end{equation}
\end{lemma}
The proof of this lemma is in appendix \ref{sec:LemmaProof}. 
Then, using this lemma the bound $\bm{\tilde \Sigma}$ on the mean
square error is 
propagated between  observation updates by solving
(\ref{eq:DtildeSig}). The updated value $\bm{\hat \Sigma}$ of the
bound $\bm{\tilde \Sigma}$ after the observation processing step
serves as the initial 
condition for $\bm{\tilde \Sigma}$.
Using a larger number of $\bm{T}$'s than $n$ may make possible a
tighter bound for $\bm{\Sigma}$. In this case there exist an infinite
number of vectors $\bm{s}$ that satisfy (\ref{eq:sdef}). One
possibility is to chose $\bm{s}$ that minimizes 
\footnote{${\rm Tr}\{\bm{\tilde \Sigma}^2\}$ is equal to the sum of
squares of the eigenvalues of $\bm{\tilde \Sigma}$. Since $\bm{\tilde
  \Sigma} \succeq \bm{0}$ this leads to a type of minimum mean square
error matrix. Another possibility is to minimize the trace of
$\bm{\tilde \Sigma}$. This also leads to linear program.}
\begin{equation}\label{eq:sopt}
\frac{{\rm d}{\rm Tr}\{\bm{\tilde \Sigma}^2\}}{{\rm d}t} = 2{\rm
  Tr}\{\bm{\tilde \Sigma} \bm{\dot{\overline \Sigma}}\} =
2\sum_i s_i {\rm Tr}\{\bm{\tilde \Sigma}   \bm{T}_i\}
\end{equation}
Minimizing (\ref{eq:sopt}) subject to (\ref{eq:sdef}) is a standard
linear programing problem.
\subsection{The Calculation of $\bm{Q}_i$ and $q_i$} \label{eq:QqCalc}
The $\bm{Q}_i$ and $q_i$, required to calculate the bound 
$\bm{\dot{\overline \Sigma}}$ are determined next. The
calculation is made 
very much easier if (\ref{eq:dif}) is composed of polynomials or
polynomial fractions. In this
case the theory of sum of squares polynomials facilitates the
calculations, at a cost of being slightly conservative. For each of
the $\bm{P}_i$ a $\bm{Q}_i$ and 
 a $q_i$ needs to be calculated such that
\begin{equation}\label{eq:SOSbound}
\begin{split}
\left[\bm{f(\tilde x) - f(\tilde x - \tilde e)}\right]^\prime \bm{P}_i\bm{ \tilde
  e} + \bm{\tilde e^\prime P}_i \left[\bm{f(\tilde x) - f(\tilde x - \tilde
      e)}\right] \\
+ {\rm Tr}\{ \bm{g^\prime(\tilde x - \tilde e) P}_i \bm{
    g(\tilde x - \tilde e)}\} -  \bm{\tilde e}^\prime\bm{ Q}_i
  \bm{\tilde e} - q_i \le 0, \; \forall \bm{\tilde e}
\end{split}
\end{equation}
The left side of (\ref{eq:SOSbound}) is an expression in 
$\bm{\tilde e}$, with $\bm{\tilde x}$ a known parameter. Any
$\bm{Q}_i$ and $q_i$ satisfying the inequality leads to a valid
bound, but with different degree of conservatism.
Assumption \ref{as:Pieq} ensures the
existence of at least one $\bm{Q}_i$ and $q_i$ for any $\bm{P}_i$. 
Among all the possible $\bm{Q}_i$ and $q_i$ that satisfy
(\ref{eq:SOSbound}), a
reasonable choice for $\bm{Q_i}$ and $q_i$ are the values that
minimize 
\begin{equation}
{\rm \bf E}\{\bm{\tilde e}^\prime\bm{Q}_i\bm{\tilde e} + q_i\} = {\rm
  Tr}\{\bm{Q}_i \bm{\Sigma}\} + q_i
\end{equation}
but $\bm{\Sigma}$ is not available. A  alternative is to
replace $\bm{\Sigma}$ with $\bm{\tilde \Sigma}$ and minimize
\begin{equation}\label{eq:SOSOpt}
{\rm Tr}\{\bm{Q}_i \bm{\tilde\Sigma}\} + q_i
\end{equation}
subject to (\ref{eq:SOSbound}).  This entails
doing the minimization in 
real time, at each integration step.  
When (\ref{eq:dif}) is polynomial or consists of polynomial fractions,
the minimization is
efficiently solved using \cite{SOSTOOLS}. In general
the calculated $\bm{Q}_i$ and $q_i$ are functions of $\bm{\tilde x}$
(as well as $\bm{\tilde \Sigma}$). The calculated $\bm{Q}_i$ and $q_i$
are then used in (\ref{eq:Trineq}). 

In addition to solving the optimization (\ref{eq:SOSOpt}), the maximum
in (\ref{eq:t_idef}) must be calculated. This, however, is a standard
linear programming problem.

\subsection{The Linear Case}
As a first demonstration of the calculation of the bound, the
algorithm is applied to the linear case. In this case the diffusion
(\ref{eq:dif}) becomes
\begin{equation}
{\rm d}\bm{x} = \bm{Ax}{\rm d}t + \bm{B}{\rm d}\bm{w}
\end{equation}
Assumptions \ref{as:ContDif} and \ref{as:Pieq} are trivially
satisfied.  For
the linear case, (\ref{eq:SOSbound}) becomes,
\begin{equation}\label{eq:SOSlinbound}
\bm{\tilde e}^\prime\bm{A}^\prime\bm{P_i\tilde e}+\bm{\tilde e}^\prime
\bm{P_i A \tilde e}
+ {\rm Tr}\{\bm{B}^\prime \bm{P_i B}\} -\bm{\tilde e}^\prime \bm{Q_i
  \tilde e} -q_i \le 0 
\end{equation}
Minimal $\bm{Q_i}$ and $q_i$ satisfying (\ref{eq:SOSlinbound})
\begin{equation}
\bm{Q_i} = \bm{A}^\prime \bm{P_i} + \bm{P}_i\bm{A}
\end{equation}
and
\begin{equation}
q_i = {\rm Tr}\left\{\bm{B}^\prime \bm{P}_i\bm{B}\right\}
\end{equation}
Substituting these into (\ref{eq:Trineq}), and rearranging,
\begin{equation}
{\rm Tr}\left\{\bm{P_i\left(\dot \Sigma -\Sigma A^\prime - A \Sigma- B
  B^\prime\right)}\right\} \le 0
\end{equation}
which, in view of (\ref{eq:Pidef}) implies
\begin{equation}
\bm{\dot \Sigma =\Sigma A^\prime + A \Sigma + B B^\prime}
\end{equation}
Hence, for the linear case, the method of the bounds employed in the
present paper in fact implies the exact equations for the propagation
of the covariance.
\section{The Observation Update Step}
At the start of the observation
update step only a bound on the mean square error is available. Using
this bound both the estimate and the mean square error bound is to be
updated. The update step is
 restricted to be linear and of the form
\begin{equation}\label{eq:update}
\bm{\hat x}_k = \bm{\tilde x}(T_k) + \bm{K}_k\left(\bm{y}(T_k) - \bm{H
  \tilde x}(T_k)\right) 
\end{equation}
where $\bm{K}_k$ is a gain matrix, which is to be determined.
A reasonable criterion for the updated estimate is one that minimizes
the worst case mean squared error given the bound  $\bm{\tilde
  \Sigma}(T_k)$ on the prior mean squared error matrix $\bm{\Sigma}(T_k)$. The
expected value 
of the mean squared error matrix  after the update (\ref{eq:update})
is given by the Joseph 
form, \cite{BarShalom}.
\begin{equation}
\bm{\hat \Sigma}_k = \left[\bm{I - K_kH\right}]\bm{\Sigma}(T_k)\left[\bm{I -
    K_kH}\right]^\prime + \bm{K}_k\bm{RK}_k^\prime
\end{equation}
Accordingly, $\bm{K}_k$ is chosen as the value that achieves the mini-max in
\begin{equation}
\min_{\bm{K}_k}\max_{\bm{\Sigma}(T_k) \preceq \bm{\tilde
    \Sigma}(T_k)} {\rm Tr}\{\bm{\hat \Sigma}_k\}
\end{equation}
The $\bm{\Sigma}(T_k)$ that achieves the inner maximum of the trace is $\bm{\tilde
  \Sigma}(T_k)$. The value of $\bm{K}_k$  that minimizes the trace of the Joseph
form is given by,  
\begin{equation}\label{eq:Kgain}
\bm{K}_k = \bm{\tilde \Sigma}(T_k)\bm{H}^\prime[\bm{R +
    H\tilde{\Sigma}}(T_k) \bm{H}^\prime]^{-1} 
\end{equation}
This then is the well known Kalman filter gain, but with the mean
squared error matrix replacing the prior covariance matrix. The
updated bound on the mean squared error matrix is 
\begin{equation}\label{eq:SigmaUpdate}
\bm{\hat \Sigma}_k = \left[\bm{I -
    K}_k\bm{H}\right]\bm{\tilde\Sigma}(T_k)\left[\bm{I - 
    K}_k\bm{H}\right]^\prime + \bm{K}_k\bm{RK}_k^\prime
\end{equation}

\section{Implementing the Algorithm}
In this section the various parts of the filtering algorithm are
collected.
\subsection{Off Line Calculations}
The $\bm{P}_i$ and $\bm{T}_i$ are independent of the observations, and
the filter state, therefore these matrices may be calculated
off-line. The algorithm for the calculation can be found in appendices
\ref{sec:PiCalc} and \ref{sec:TiCalc}. Also $\bm{L}$ is calculated
at this stage (\ref{eq:Lcalc}). 
\subsection{Observation Update}
The observation update very much follows the observation update of the
conventional extended Kalman filter, except the prior covariance
matrix is replaced by the mean square error matrix bound $\bm{\tilde
  \Sigma}$. The state is updated as in (\ref{eq:update}), the gain is
calculated in (\ref{eq:Kgain}), and the mean square error matrix is
updated in (\ref{eq:SigmaUpdate}).
\subsection{Time Update}
Between observations the state is propagated according
to (\ref{eq:Tupdt}). The mean square error matrix is propagated by
solving (\ref{eq:DtildeSig}). In order to be able to carry this out
the matrix function 
$\bm{\dot {\overline \Sigma}}(\bm{\tilde \Sigma},t)$ is required. This
function is calculated at each integration step of
(\ref{eq:DtildeSig}). The $\bm{\hat\Sigma}_k$ calculated in
(\ref{eq:SigmaUpdate}) serves as the initial condition for
(\ref{eq:DtildeSig}) after the $k$th data processing step.
\begin{enumerate}
\item The values $\bm{\tilde x}$,  $\bm{Q}_i$ and $q_i$ are
  calculated, by 
minimizing the expression (\ref{eq:SOSOpt}) subject to
(\ref{eq:SOSbound}). When the process (\ref{eq:dif}) consists of
polynomials or rational functions of polynomials this step is most
efficiently carried out using \cite{SOSTOOLS}. An example of this
calculations is presented in section \ref{sec:example}.
\item The 
linear program (\ref{eq:t_idef}), subject to (\ref{eq:SetS}), with
$\bm{\Sigma}=\bm{\tilde\Sigma}$ and the minimization (\ref{eq:sopt}) subject to
(\ref{eq:sdef})  are solved to calculate $\bm{s}$ and finally  
$\bm{\dot{\overline \Sigma}}$ is obtained by substituting $\bm{s}$ in 
(\ref{eq:bardotSig}). 
\end{enumerate}
\section{An Example}\label{sec:example}
The example in the paper is a second order system with a limit cycle,
\begin{equation}\label{eq:exampl}
{\rm d}\bm{x} = \bm{A}(\bm{x})\bm{x} +\bm{G}{\rm d}\bm{w}
\end{equation}
where
\begin{equation}
\bm{A}(\bm{x}) = \bm{A}_u/m(\bm{x}) + \bm{A}_s
\end{equation}
\begin{equation}
\bm{A}_u = \left[
\begin{array}{cc}
1 & 1\\
-1 & 1
\end{array}\right], \;
\bm{A}_s = \left[
\begin{array}{cc}
-1 & 1\\
-1 & -1
\end{array}\right]
\end{equation}
\begin{equation}
m(\bm{x}) = (1 + \bm{x}^\prime \bm{x})/25
\end{equation}
\begin{equation}
\bm{G} = \left[
\begin{array}{cc}
1/5 & 0\\
0  & 1/5
\end{array}\right]
\end{equation}
The observation matrix is
\begin{equation}
\bm{H} = \left[
\begin{array}{lr}
1 & 0
\end{array}\right]
\end{equation}
and $y = \bm{Hx} + v_k$, is observed every $0.2$s. The standard
deviation of $v_k$ is $0.01$.

The initial step in calculating the bound is
substituting into (\ref{eq:SOSbound}) yields
\begin{equation}\label{eq:SOSbndExample}
\begin{split}
\left[\bm{A}(\bm{\tilde x})\bm{\tilde x} - \bm{A}(\bm{\tilde x -
    \tilde e})(\bm{\tilde x - \tilde
    e})\right]^\prime\bm{P}_i\bm{\tilde e} \\+
\bm{\tilde e}^\prime\bm{P}_i\left[\bm{A}(\bm{\tilde x})\bm{\tilde x} -
  \bm{A}(\bm{\tilde x - \tilde e})(\bm{\tilde x - \tilde e})\right]\\
+{\rm Tr}\{\bm{G^\prime} \bm{P}_i \bm{G}\} - \bm{\tilde e}^\prime
\bm{Q}_i \bm{\tilde e} - q_i
\end{split}
\end{equation}
This form is not suitable for use in SOSTOOLS because of the form
$m(\bm{\tilde x - \tilde e})$ in the denominator. Since 
$m(\bm{\tilde  x - \tilde e})> 0$ multiplying
(\ref{eq:SOSbndExample}) by $m(\bm{\tilde x - \tilde e})$ does not
affect the sign of the form, so that (\ref{eq:SOSbndExample}) may be
replaced by
\begin{equation}\label{eq:SOSboundExample2}
\begin{split}
m(\bm{\tilde x - \tilde e})\left[\bm{A}(\bm{\tilde x})\bm{\tilde x} -
  \bm{A}(\bm{\tilde x - \tilde e})(\bm{\tilde x - \tilde
    e})\right]^\prime\bm{P}_i\bm{\tilde e} \\+
m(\bm{\tilde x - \tilde e})\bm{\tilde
  e}^\prime\bm{P}_i\left[\bm{A}(\bm{\tilde x})\bm{\tilde x} - 
  \bm{A}(\bm{\tilde x - \tilde e})(\bm{\tilde x - \tilde e})\right]\\
+m(\bm{\tilde x - \tilde e}){\rm Tr}\{\bm{G^\prime} \bm{P}_i \bm{G}\}
- m(\bm{\tilde x - \tilde e})\left(\bm{\tilde e}^\prime 
\bm{Q}_i \bm{\tilde e} - q_i \right)
\end{split}
\end{equation}
Then the input to SOSTOOLS consists of the minimization of
(\ref{eq:SOSOpt})  subject to the constraint
(\ref{eq:SOSboundExample2}).

\begin{figure}
\centering
\includegraphics[width = 0.8\columnwidth]{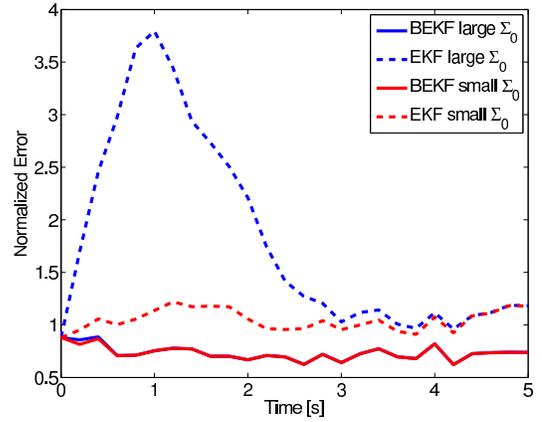}
\caption{Average normalized mean square error 
$\sqrt{\bm{\tilde e}^\prime\bm{ \tilde \Sigma}^{-1}\bm{\tilde e}/2}$}
\label{fig:MSE}
\end{figure}
In Fig.  \ref{fig:MSE} is shown the average normalized mean square error,
$\sqrt{\bm{\tilde e}^\prime\bm{ \tilde \Sigma}^{-1}\bm{\tilde e}/2}$,
calculated at each integration step, for the initial
conditions 
\begin{equation}\label{eq:case1}
\bm{\Sigma}_0 = \left[
\begin{array}{cc}
0.5 & 0\\
0 & 0.5
\end{array}\right], \; \bm{\mu} = \left[
\begin{array}{c}
8\\
0
\end{array}\right]
\end{equation}
\begin{equation}\label{eq:case2}
\bm{\Sigma}_0 = \left[
\begin{array}{cc}
0.01 & 0\\
0 & 0.01
\end{array}\right], \; \bm{\mu} = \left[
\begin{array}{c}
8\\
0
\end{array}\right]
\end{equation}
Two sets of lines are shown in Fig \ref{fig:MSE}: the full lines (blue
and red) shows the
normalized error for the BEKF filter proposed in the
present paper; and the dashed lines (blue and red) the normalized
error for the 
classical extended Kalman filter (here the error is normalized using
the covariance matrix calculated by the extended Kalman filter). It is
seen that for the proposed filter the normalized error is for all
cases about $0.7$, implying a certain amount of
conservativeness. For the extended Kalman filter, when the initial
uncertainty is small red line, (\ref{eq:case2}), the normalized error
is near 
one, implying that the filter is operating well. For case
(\ref{eq:case1}), blue dashed line, the conventional extended Kalman
filter severely underestimates the filter error, implying considerably
too low a value for the filter gain. Only after about 3
seconds is does the normalized mean square error settle near 1.
 Further increasing the initial
uncertainty for the extended Kalman filter causes it to loose all
connection between the calculated covariance matrix and the actual
estimate error. Note that this situation may be remedied by increasing
the process noise in the extended Kalman filter beyond its actual value.

\section{Summary and Conclusion}
A new form of the extended Kalman filter was presented. The main
distinguishing feature of the new filter is the computation of a bound
for the mean
square error matrix for the estimate error. This matrix serves as both
a bound on the actual mean square error and is used in the calculation
of the Kalman gain. In contrast, the connection between the covariance
matrix, computed by linearizing the dynamics about the state estimate
in the conventional extended Kalman filter
is at best very approximate. As a consequence, the extended Kalman
has to be ``tuned'' by increasing the process noise amplitude, to
ensure some sort of connection between the computed covariance and the
actual mean square error.

The emphasis has been on deriving a mean square bound. Tighter bounds
may be achievable by judicious choice of the matrices $\bm{P}_i$,
and $\bm{U}_i$. Note that when the system is linear the proposed
filter reduces to the conventional linear Kalman filter.
\appendix
\subsection{Calculating $\bm{T}_i$}\label{sec:TiCalc}
The purpose of this section is to calculate the $\bm{T}_i$ from the
$\bm{U}_i$. As noted in section \ref{sec:dotSigma} the bound calculated in the paper
becomes tighter the better the cone generated by the set of $\bm{U}_i$
approximate ${\cal M_+}$, the cone of positive definite
matrices. There is at present no explicit algorithm for selecting the
$\bm{U}_i$, except for choosing the $\bm{U}_i$ to be rank one positive
semi-definite matrices. Note that {\em any} choice of a set of positive
semi-definite matrices, $\bm{U}_i$ results in a valid mean square
error bound. but some choices for $\bm{U}_i$ may result in a tighter bound.

When there are $n$ $\bm{U}_i$'s, to
each $\bm{U}_i$ there corresponds a 
$\bm{T}_i$ and is defined as follows
\begin{enumerate}
\item $\bm{T}_i$ is orthogonal  to all the $\bm{U}_k,\; k \not = i$
that is ${\rm Tr}\{\bm{T}_i\bm{U}_k\} = 0$.
\item in addition, ${\rm Tr}\{\bm{T}_i\bm{U}_i\} =1 > 0$
\end{enumerate}
These two conditions define $n$ linear equations in $n$ unknowns,
whose solution is straightforward.
When there are more than $n$ $\bm{U}_i$'s then each subset of $n$
$\bm{U}_i$'s yields a set of $\bm{T}_i$'s. From the union of all these
$\bm{T}_i$'s corresponding to all possible subsets of $n$
$\bm{U}_i$'s, a set of $\bm{T}_i$'s is chosen as the generators of
${\cal F}^*$
\begin{enumerate}
\item If ${\rm Tr}\{\bm{T}_i\bm{U}_j\} > 0$ for some $j$ then
  $\bm{T}_i$ is excluded.
\item Any $\bm{T}_i$ not in the set is expressible as a positive
  linear combination of the $\bm{T}_i$'s in the set.
\item Any $\bm{T}_i$ in the set is not expressible as a positive
  linear combination of the remaining $\bm{T}_i$'s in the set.
\end{enumerate}
\subsection{Calculating the $\bm{P}_i$}\label{sec:PiCalc}
This section presents an algorithm for calculating a minimal set of
$\bm{P}_i$'s. 
The main  property of the $N\times N$ symmetric matrices $\bm{P}_i$ is
that given {\em any} $N\times N$ symmetric 
matrix $\bm{Z}$ there exists $i$ such that ${\rm Tr}\{\bm{P}_i\bm{Z}\}>0$.
Let $\bm{P}_1$ be arbitrary, subject to
\begin{equation}
{\rm Tr}\{\bm{P}_1 \bm{P}_1\} = 1
\end{equation}
There are $n+1$ $\bm{P}$'s, and the first $n$  are defined recursively.
\begin{equation}
\bm{P}_{k+1} = \alpha_k \sum_{i=1}^k\bm{P}_i + \beta_k
\bm{S}_k,\;k=1,2,\ldots, n-1
\end{equation}
The $\bm{S}_k$ is chosen to satisfy
\begin{equation}
{\rm Tr}\{\bm{P}_i \bm{S}_k\} = 0, \;\text{for } i=1,2\ldots k
\end{equation}
and
\begin{equation}
{\rm Tr}\{\bm{S}_k \bm{S}_k\} = 1
\end{equation}
Note that some of the entries of $\bm{S}_k$ are arbitrary, since the
number of entries in $\bm{S}_k$ is greater than the number of equations.
The constant $\alpha_k$ is chosen so that
\begin{equation}\label{eq:PiPj_1n}
{\rm Tr}\{\bm{P}_{k+1} \bm{P}_i\} = -1/n,\;\text{ for } i = 1,2,\ldots
k
\end{equation}
Giving $\alpha_k$ the value
\begin{equation}
\alpha_k = -\frac{1}{n-(k-1)}
\end{equation}
ensures this. The constant $\beta_k$ is chosen so that
\begin{equation}\label{eq:PiPi_1}
{\rm Tr}\{\bm{P}_{k+1} \bm{P}_{k+1}\} = 1
\end{equation}
Choosing
\begin{equation}
\beta_k = \sqrt{1-\frac{k}{n(n-(k-1))}}
\end{equation}
ensures this. Note that the expression under the radical is always
positive. Finally the last $\bm{P}_{n+1}$ is defined
\begin{equation}\label{eq:Pn_1def}
\bm{P}_{n+1} = -\sum_{i=1}^n \bm{P}_i
\end{equation}
Using (\ref{eq:PiPi_1}) and (\ref{eq:PiPj_1n}) and the definition of
$\bm{P}_{n+1}$ in (\ref{eq:Pn_1def}) yields
\begin{equation}
{\rm Tr}\{\bm{P}_{n+1} \bm{P}_{n+1}\} = 1
\end{equation}
A similar calculation shows that
\begin{equation}
{\rm Tr}\{\bm{P}_{n+1} \bm{P}_i\} = -1/n,\; i=1,2,\ldots, n
\end{equation}
From their construction the first $n$ $\bm{P}_i$ span the space 
of $N \times N$ symmetric matrices. Accordingly, if $\bm{Z}$ is any
symmetric $N \times N$ matrix there exist $z_i$,
\begin{equation}\label{eq:Pexp}
\bm{Z} = \sum_{i=1}^n z_i \bm{P}_i
\end{equation}
Solving (\ref{eq:Pn_1def}) for any $\bm{P}_i$ and substituting into
(\ref{eq:Pexp}) shows that in fact any $n$ $\bm{P}_i$'s span the space
$N \times N$ symmetric matrices. 
Rearranging (\ref{eq:Pn_1def})
\begin{equation}
\sum_{i=1}^{n+1} \bm{P}_i = \bm{0}
\end{equation}
so that
\begin{equation}\label{eq:ZPi}
{\rm Tr}\left\{\bm{Z} \sum_{i=1}^{n+1} \bm{P}_i\right\} =0
\end{equation}
Since the $\bm{P}_i$'s span the space of symmetric matrices, if
$\bm{Z}\not = \bm{0}$, not all the terms in (\ref{eq:ZPi}) are
zero. Therefor there exists an $i$ such that
\begin{equation}\label{eq:Zpos}
{\rm Tr}\{ \bm{Z} \bm{P}_i\} > 0
\end{equation}
\subsection{Proof of Lemma \ref{lem:comp}}\label{sec:LemmaProof}
Suppose that the Lemma is false. Then there exists $t_1$, and $t_2$
such that
\begin{equation}
\bm{ \Sigma}(t) \succ \bm{\tilde\Sigma}(t),\; t_1 < t \le t_2
\end{equation}
but
\begin{equation}
\bm{ \Sigma}(t_1) = \bm{\tilde\Sigma}(t_1)
\end{equation}
From the mean value theorem,
\begin{equation}
\bm{\tilde \Sigma}(t_2) - \bm{ \Sigma}(t_2)
=\frac{{\rm d}}{{\rm d}t}\left( \bm{\tilde \Sigma}(t_a) - \bm{
    \Sigma}(t_a)\right)(t_2 -t_1)
\end{equation}
where $t_1 \le t_a \le t_2$, implying that 
\begin{equation}
\frac{{\rm d}}{{\rm d}t}\left( \bm{\tilde \Sigma}(t_a) - \bm{
    \Sigma}(t_a)\right) \succ \bm{0}
\end{equation}
which contradicts (\ref{eq:DtildeSig}) and (\ref{eq:dotbarSigdef}).
\bibliography{IEEEabrv,Bounds}

\begin{thebibliography}{10}
\providecommand{\url}[1]{#1}
\csname url@rmstyle\endcsname
\providecommand{\newblock}{\relax}
\providecommand{\bibinfo}[2]{#2}
\providecommand\BIBentrySTDinterwordspacing{\spaceskip=0pt\relax}
\providecommand\BIBentryALTinterwordstretchfactor{4}
\providecommand\BIBentryALTinterwordspacing{\spaceskip=\fontdimen2\font plus
\BIBentryALTinterwordstretchfactor\fontdimen3\font minus
  \fontdimen4\font\relax}
\providecommand\BIBforeignlanguage[2]{{%
\expandafter\ifx\csname l@#1\endcsname\relax
\typeout{** WARNING: IEEEtran.bst: No hyphenation pattern has been}%
\typeout{** loaded for the language `#1'. Using the pattern for}%
\typeout{** the default language instead.}%
\else
\language=\csname l@#1\endcsname
\fi
#2}}

\bibitem{PapaPrajna}
A.~Papachristodoulou and S.~Prajna, ``A tutorial on sum of squares techniques
  for system analysis,'' in \emph{American Control Conference}, 2005.

\bibitem{ParriloLall}
P.~A. Parrilo and S.~Lall, ``Semidefinite programming relaxations and algebraic
  optimization in control,'' \emph{European Journal of Control}, vol.~9, no.
  2--3, pp. 307--321, 2003.

\bibitem{Chesi}
G.~Chesi, ``{LMI} techniques for optimization over polynomials in control: a
  survey,'' \emph{{IEEE} Transactions on Automatic Control}, vol. AC-55,
  no.~11, pp. 2500--2510, 2010.

\bibitem{SOSTOOLS}
\BIBentryALTinterwordspacing
S.~Prajna, A.~Papachristodoulo, P.~Seiler, and P.~A. Parrilo. (2004)
  {SOSTOOLS}. [Online]. Available: \url{http://www.cds.caltech.edu/sostools/}
\BIBentrySTDinterwordspacing

\bibitem{SeDuMi}
\BIBentryALTinterwordspacing
J.~F. Strum. (2010) {SeDuMi}. [Online]. Available:
  \url{http://sedumi.ie.lehigh.edu/}
\BIBentrySTDinterwordspacing

\bibitem{Sorenson}
H.~W. Sorenson, ``On the development of practical nonlinear filters,''
  \emph{Information Sciences}, vol.~7, pp. 253--270, 1974.

\bibitem{Jazwinsky}
A.~H. Jazwinsky, \emph{Stochastic Processes and Filtering Theory}.\hskip 1em
  plus 0.5em minus 0.4em\relax New York: Academic Press, 1970.

\bibitem{KumarChakra}
M.~Kumar and S.~Chakravorty, ``Nonlinear filter based on the {F}okker--{P}lanck
  equation,'' \emph{Journal of Guidance, Control and Dynamics}, vol.~35, no.~1,
  pp. 68--79, January--February 2012.

\bibitem{GilmanRhodes1}
A.~S. Gilman and I.~A. Rhodes, ``Cone-bounded nonlinearities and mean-square
  bounds --- quadratic regulation bounds,'' \emph{{IEEE} Transactions on
  Automatic Control}, vol. AC-21, no.~4, pp. 472--483, August 1976.

\bibitem{GilmanRhodes2}
------, ``Cone-bounded nonlinearities and mean square bounds---estimation lower
  bound,'' \emph{{IEEE} Transactions on Automatic Control}, vol. AC-20, no.~5,
  pp. 632--641, 1975.

\bibitem{GilmanRhodes3}
------, ``Cone-bounded nonlinearities and mean square bounds---estimation upper
  bound,'' \emph{{IEEE} Transactions on Automatic Control}, vol. AC-18, no.~3,
  pp. 260--265, 1973.

\bibitem{BarShalom}
Y.~BarShalom, X.~R. Li, and T.~Kirubajan, \emph{Estimation with Applications to
  Tracking and Navigation: Theory Algorithm and Software}.\hskip 1em plus 0.5em
  minus 0.4em\relax John Wiley and Sons Inc., 2001.

\bibitem{Lohmiller}
W.~Lohmiller and J.~E. Slotine, ``On contraction analysis for non-linear
  systems,'' \emph{Automatica}, vol.~34, no.~6, pp. 683--696, 1998.

\bibitem{Pham}
Q.~Pham, N.~Tabareau, and J.~E. Slotine, ``A contraction theory approach to
  stochastic incremental stability,'' \emph{{IEEE} Transactions on Automatic
  Control}, vol. AC-54, no.~4, 2009.

\bibitem{Boyd}
S.~Boyd and L.~Vandenberghe, \emph{Convex Optimization}.\hskip 1em plus 0.5em
  minus 0.4em\relax Cambridge University Press, 2004.

\bibitem{Khalil}
H.~K. Khalil, \emph{Nonlinear Systems}, 3rd~ed.\hskip 1em plus 0.5em minus
  0.4em\relax Prentice Hall, 2002.

\end{thebibliography}
\end{document}